\documentclass[12pt]{article}
\title{Very Strong Disorder for the Parabolic Anderson model in low dimensions}
\author{Pierre Bertin}

\usepackage{amssymb}
\usepackage[mathcal]{eucal}

\usepackage{enumerate}
\usepackage{palatino}

\usepackage{pgf,tikz}
\usetikzlibrary{arrows}

\definecolor{ffqqqq}{rgb}{1,0,0}
\definecolor{qqqqff}{rgb}{0,0,1}
\definecolor{cqcqcq}{rgb}{0.75,0.75,0.75}
\definecolor{qqffqq}{rgb}{0,1,0}

\usepackage{amssymb}
\usepackage{bbm}

\usepackage{amsthm}

\newtheorem{theo}{Theorem}

\newcommand{\R}{\mathbb{R}}

\newcommand{\Z}{\mathbb{Z}}

\renewcommand{\Z}{\mathbb{Z}}

\def\ind{{\mathbbm{1}}_}

\def\vep{\varepsilon}

\def\m{{\mathcal M}}

\def\g{{\mathcal G}}

\def\z{{\mathcal Z}}

\def\be{\begin{equation}}
\def\ee{\end{equation}}
\def\ba{\begin{eqnarray*}}
\def\ea{\end{eqnarray*}}

\def\wh{\widehat}
\def\wt{\widetilde}

\begin{document}

\begin{center} \begin{LARGE}
Very Strong Disorder for the Parabolic Anderson model in low dimensions
\end{LARGE}\end{center}
\vspace{1cm}

\begin{large} Pierre Bertin \end{large} \\
\begin{small} Universit\'e Diderot - Paris 7, LPMA Case 7012 Paris Cedex 13 France \\ Ecole Normale Sup\'erieure de Paris, DMA 45 rue d'Ulm 75005 Paris France \\
pierre.bertin@ens.fr \end{small}

\section*{Abstract}
We study the free energy of the Parabolic Anderson Model, a time-continuous model of directed polymers in random environment. We prove that in dimension 1 and 2, the free energy is always negative, meaning that very strong disorder always holds. The result for discrete polymers in dimension two, as well as better bounds on the free energy on dimension 1, were first obtained by Hubert Lacoin in \cite{lac}, and the goal of this paper is to adapt his proof to the Anderson Parabolic Model.

\vspace{0.5cm}
Keywords : Parabolic Anderson model; Brownian Directed Polymer; Free energy; Strong disorder; Coarse graining.

\section{Model and known results}

\subsection{The Parabolic Anderson Model}

The general model we consider in this paper is defined in terms of a random trajectory in a random environment :
\begin{itemize}
	\item The random trajectory : let $((X_t)_{t\geq 0}, P^\kappa)$ be a nearest-neighbour random walk on the lattice $\Z^d$ with rate of jump $\kappa$ starting from 0. Specifically, we let $(\m,\mathcal{F})$ be the path space of c\`adl\`ag trajectories from $\R_+$ to $\Z^d$, and the process $(X_t,P^\kappa)$ has for infinitesimal generator $\kappa\bigtriangleup^d$, where $\bigtriangleup^d$ is the discrete Laplacian. Sometimes we will need to start the random walk from a point other than 0, let $((X_t),P_x^\kappa)$ be the walk starting from the point $x$. 

	\item The random environment : let $((B_t^x(\omega))_{t\geq 0}, x\in \Z^d)$ be a family of independent Brownian motions defined on a probability space $(\Omega, \mathcal{G},Q)$, one for each vertex of $\Z^d$. For $t>0$, it is natural and convenient to introduce the sub-$\sigma$-field 
\[ \mathcal{G}_t = \sigma\{ B_s^x,s \in [0,t], x\in \Z^d \ \}. \]

	\item The Gibbs measure : for a time horizon $t$ we define the following Hamiltonian :
\[ H_t = H_t (X) = \int_0^t dB_s^{X_s}= \sum_{x\in\Z^d}\int_0^t \ind{\{X_s=x\} } dB^x_s. \]
For any $t>0$, we define a probability measure $\mu_t^{\kappa,\beta}$ on the path space $(\Omega,\mathcal{F})$
\[ \mu^{\kappa,\beta}_t(dX)= \frac{\exp (\beta H_t(X)-t\beta^2/2)}{Z^{\kappa,\beta}_t}P^\kappa(dX), \]
where $\beta \geq 0$ is the inverse temperature and 
\[ Z^{\kappa,\beta}_t = P^\kappa\left[\exp(\beta H_t(X)-t\beta^2/2)) \right]. \]
Since the polymer measure is parametrized by the environment $\omega$, it is random. 
\end{itemize}

This model has been studied in \cite{camo} with the point of view of intermittency, and in \cite{gaho} when interacting with a particle system. It has also been studied as an example of polymer model in \cite{cativi}, \cite{cahu}. In \cite{mooc} Moriarty and O'Connell studied a perfectly assymetrical version of the 1-dimensional model in which they could prove exact formulas.

In this paper, we investigate the asymptotic behavior of the partition function $Z^{\kappa,\beta}_t$ when $t\to\infty$.

Let us finish the definition of the model with some remarks on scales. The process $(X_t,P^\kappa)$ has the same law than $(X_{\kappa t},P^1)$, and if $B_t$ is a Brownian motion, $B_{\kappa t}$ has the same law than $\sqrt{\kappa}B_t$. Therefore, $Z^{\kappa,\beta}_t$ has the same law than $Z^{1,\frac{\beta}{\sqrt{\kappa}}}_{\kappa t} $, and we only need to study the case of $\kappa=1$, and the rest is obtain by rescaling.

We will denote by $P$, $\mu^\beta_t$, $Z^\beta_t$,... the quantities $P^\kappa$, $\mu^{\kappa,\beta}_t$, $Z^{\kappa,\beta}_t$,... with $\kappa=1$. In most case, when the value of $\beta$ is not ambiguous, we will only write $\mu_t$ and $Z_t$.

\subsection{The partition function}

Let us fix a $\beta\geq 0$ and begin by stating that the partition function $Z_t$ is a martingale on $(\mathcal{M}, \mathcal{G},Q)$. For any fixed path $X$, the process $\{H_t (X)\}_{t\geq 0}$ is himself a Brownian motion, and $\{ \exp(\beta H_t - t\beta^2/2)\}_{t\geq 0}$ is its exponential martingale. Therefore, the partition function $Z_t$
is a mean-one, continuous, positive martingale on $(\mathcal{M}, \mathcal{G},Q)$ with respect to the filtration $(\mathcal{G}_t)_{t\geq 0}$. In particular, the following limit exists $Q$-a.s. :
\[ Z_\infty = \lim_{t\to \infty} Z_t \]
Since $\exp (\beta H_t)>0$, the event $\{Z_\infty =0\}$ is measurable with respect to the tail $\sigma$-field 
\[\bigcap_{t\geq1} \sigma\{ B^x_s, s>t, x\in \Z^d \} \]
and therefore by Kolmogorov's $0-1$ law, 
\[ Q\{Z_\infty =0\} = 0 \textrm{ or } 1.\]
In the case $Z_\infty>0$, we say that \emph{weak disorder} holds, and in the case $Z_\infty=0$ that \emph{strong disorder} holds.

\subsection{The phase transition}

\begin{theo} There exists a critical value $\beta_c=\beta_c(d)\in[0,\infty)$ such that :
\begin{itemize} 
\item if $\beta<\beta_c$ then $Z_\infty >0$ a.s.
\item if $\beta>\beta_c$ then $Z_\infty =0$ a.s.
\end{itemize}
Moreover, if $d\geq 3$ then $\beta_c(d)>0$. \end{theo}

A lot of work has been made to study these two different phases for random walk in random environment models, and it has strong links with the overlap of two independent trajectories in the same environment $\omega$. If we call $J_t$ the expected overlap of two independent replicas :
\[J_t=J_t(\beta)=\mu_t^{\beta \otimes 2}\left[\int_0^t \ind{\{ X_s=\widehat{X}_s\}}ds\right],\]
then,  roughly speaking :
\[ \ln Z_t^\beta \approx Q[\ln Z_t^\beta]=-\int_0^\beta bQ[J_t(b)]db.\]
For precise results on the localization of the polymer, favorite end-point and favorite trajectory, I invite you to consult \cite{cahu} and \cite{cocr}.

Empirically, in the strong disoder case, the polymer will be compelled to follow a most attractive path, while in the weak disorder phase, the probability will be more fairly distributed amongst all possible trajectories. This is why this phase transition is sometimes called the transition from delocalized phase to localized phase.

It is also interesting to know when \emph{very strong disorder} holds, meaning that $Z_t$ converges to 0 exponentially fast. It is a well-known fact in the PAM that the limit 
\[\Psi(\kappa,\beta)= \lim_{t\to\infty}\frac{1}{t}\ln Z^{\kappa,\beta}_t\]
exists a.s. and in $\mathbf{L}^p$ for $p\in[1,\infty)$. It is called the \emph{free energy} of the model. For more results on the free energy, consult \cite{crmosh}, and read \cite{crgamo} for the large deviations point of view.

\subsection{Main theorem}

In this paper we prove that very strong disorder always holds in dimension 1 and 2.

\begin{theo}
In the Parabolic Anderson Model for dimension 1 and 2, for all $\beta >0$, there exists $c>0$ such that
\[ Z_t = O(e^{-ct}) \ \ Q\mathrm{-a.s.} \]
Moreover, we have the following upper bounds for small $\beta$ ($c_1,c_2$ are positive constants):
\[\Psi(1,\beta)\leq \left\{ \begin{array}{ll} -c_1\beta ^4 & \textrm{ for }\ d=1 \\ -\exp\left(\frac{-c_2}{\beta^4}\right) & \textrm{ for }\ d=2\end{array} \right.\]
\end{theo}

In other terms, as $t\to\infty$, the partition function grows exponentially slower than its expectation.

\emph{Remark 1 :} In higher dimensions, as a corollary of theorem 1, the free energy $\Psi(1,\beta)$ is constant equal to 0 on the (non trivial) interval $[0,\beta_c(d))$ (at least).

\emph{Remark 2 :} A lot of work has also be done to study the asymptotics of $\psi(1,\beta)$ for $\beta$ big. It is well known in the Parabolic Anderson model (see \cite{cativi} \cite{cako} \cite{cocr} \cite{crmosh}) that when $\beta \to \infty$
\[\Psi(1,\beta)= \frac{\beta^2}{2} - \frac{\alpha^2 \beta^2}{4\log(\beta^2)}+ o \left(\frac{\beta^2}{\log(\beta^2)}\right),\]
where $\alpha$ is a constant depending on the dimension.

\emph{Remark 3 :} In \cite{mooc}, O'Connell and Moriarty studied a perfectly assymetrical version of the 1-dimensional model : the random walk $X_t$ always jumps in the same direction. They were able, thanks to reversibility properties of this new model, to prove an exact formula for the free energy. The asymptotics of this free energy when $\beta\to 0$ are $-\frac{\beta^4}{24}+O(\beta^6)$, which is the order found in our theorem. Moreno in \cite{mo} also showed that an assymetric Parabolic Anderson Model with very high drift had the same behavior than O'Connell and Moriarty's perfectly assymetric model.

The goal of this paper is to adapt the proof of the same results for directed polymers in random environment for discrete times, which can be found in \cite{lac} to the Parabolic Anderson Model.

The same proof can be adapted for other directed polymer models : a Brownian directed polymer in a Poisson environment as seen in \cite{coyo}, the Linear Stochastic Evolutions model introduced by Yoshida in \cite{yo1,yo2}. Similar techniques are also used for non-directed models, Lacoin shows in \cite{lac2} that the number of self avoiding walks of length $n$ in the infinite cluster of the supercritical percolation in $\Z^2$ is exponentially smaller than its expectation.

\section{Proof in dimension 1}

Before starting, we sketch the proof and how it is decomposed in different steps. These steps are the same as in the section 3 of \cite{lac}:
\begin{enumerate}
\item We reduce the problem to the exponential decay of $Q[Z_t^\theta]$ for some $\theta \in (0,1)$.
\item We use a decomposition of $Z_t$ by splitting it into different contributions that corresponds to trajectories that stay in a large corridor.
\item To estimate the fractional moment terms appearing in the decomposition, we will change the measure $Q$ around the path.
\item We use some basic properties of a random walk in $\Z$ to compute the expectation with the new measure.
\end{enumerate}

\subsection{The fractional moment method}

To show the exponential decay of $Z_t$, we will use the trick known as the fractional moment. We want to show that
\[ \lim_{t \to \infty} \frac{1}{t} Q[ \log Z_t] \leq -c < 0  \]
but it is not easy to handle the expectation of a log. So we pick $\theta \in (0,1)$, and by Jensen inequality :
\begin{equation} \label{fract} 
Q[ \log Z_t] = \frac{1}{\theta} Q[ \log Z_t^\theta ] \leq \frac{1}{\theta} \log Q[Z_t^\theta] \end{equation}
We are left with showing that $Q[Z_t^\theta] =\mathcal{O}(e^{-c't})$, which is easier. In the rest of the proof, we will take $\theta=\frac{1}{2}$.

\subsection{Decomposition of $Z_t$}

We fix $n$ to be a very large number. The meaning of very large will be seen at the end of the proof, . The number $n$ will be used in the sequel as a scaling factor. We also take $n$ such that $\sqrt{n}$ is an integer. In the following, when we define $n$, the notation $n= f(\beta)$ will mean that $n$ is the smallest squared integer bigger than $f(\beta)$. If $f(\beta)\geq 5$, then $f(\beta)\leq n <2f(\beta)$.

For $z \in \mathbb{Z}$, we define $I_z = [z\sqrt{n},(z+1)\sqrt{n}) $ so that the $I_z$ are disjoint and cover $\Z$. For $T=mn$, we sort the paths $(X_s)_{s\in[0,T]}$ of the random walk by the position of the points $(X_{in})_{i=1}^m$. More precisely, we look at the $m$-uplet $(z_1,\dots,z_m)$ such that $X_{in} \in I_{z_i}$ for $i=1,\dots,m$. 

Now we can decompose $Z_t$ :
\begin{equation} \label{decomp}
Z_t = \sum_{(z_1,\dots, z_m)\in \Z^m} \overline{Z}_{(z_1,\dots,z_m)} \end{equation}
where
\[\overline{Z}_{(z_1,\dots,z_m)} = P\left[\exp(\beta H_T(X) - T\beta^2/2) \ind{\{ X_{in} \in I_{z_i}, \ i=1\dots m \} }\right]. \]

In the sequel, for shorter notations, we will write $\mathcal{Z}=(z_1,\dots,z_m)$ and use the convention that $z_0=0$.

\begin{center}
\begin{tikzpicture}[line cap=round,line join=round,>=triangle 45,x=2.0cm,y=0.8cm]
\draw [color=cqcqcq,dash pattern=on 1pt off 1pt, xstep=2.0cm,ystep=0.8cm] (0,-1.46) grid (5.93,4.68);
\clip(-0.6,-1.46) rectangle (5.93,4.68);
\draw [line width=2pt,color=ffqqqq] (1,-1)-- (1,0);
\draw [line width=2pt,color=ffqqqq] (2,3)-- (2,2);
\draw [line width=2pt,color=ffqqqq] (3,1)-- (3,0);
\draw [line width=2pt,color=ffqqqq] (4,4)-- (4,3);
\draw [line width=2pt,color=ffqqqq] (5,0)-- (5,1);
\draw (0,0)-- (0.26,0.85)-- (0.8,-0.9)-- (1.29,0.64)-- (1.74,1.04)-- (2.1,2.78)-- (2.33,2.06)-- (2.59,0.94)-- (3.11,0.35)-- (3.4,2.69)-- (3.66,2.21)-- (4.08,3.7)-- (4.3,3.14)-- (5,0.6);
\draw [->] (0,0) -- (5.71,0);
\draw [->] (0.01,-1.58) -- (0.01,4.45);
\begin{scriptsize}
\draw[color=black] (1,0.22) node {$n$};
\draw[color=black] (3.02,-0.3) node {$3n$};
\draw[color=black] (5.05,-0.3) node {$5n$};
\draw[color=black] (-0.3,0.08) node {$0$};
\draw[color=black] (-0.3,1.06) node {$\sqrt{n}$};
\draw[color=black] (-0.3,2.05) node {$2\sqrt{n}$};
\draw[color=black] (-0.3,3.11) node {$3\sqrt{n}$};
\draw[color=black](-0.3,-0.93) node {$-\sqrt{n}$};
\draw[color=black] (2.07,-0.3) node {$2n$};
\draw[color=black] (4.07,-0.3) node {$4n$};
\draw[color=black] (5.87,0.15) node {$\R_+$};
\draw[color=black] (-0.3,4.14) node {$\Z$};
\end{scriptsize}
\end{tikzpicture}

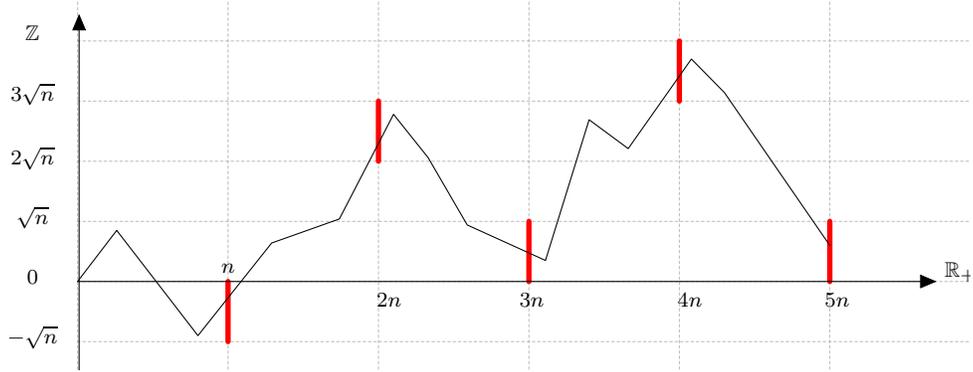
a trajectory for $\overline{Z}_{-1,2,0,3,0}$
\end{center}

Now we use our $\theta$ and the inequality $(\sum a_i)^\theta \leq \sum a_i^\theta$ (which holds for any countable collection of positive real numbers) to get :
\begin{equation} \label{decomp2} 
Q[Z_t^\theta] \leq \sum_{\z \in \mathbb{Z}^m} Q[\overline{Z}_\mathcal{Z}^\theta ] .\end{equation}

\subsection{The change of measure}

We fix one $\mathcal{Z}=(z_1,\dots,z_m)$. We define $m$ blocs $(J_k)_{k=0}^{m-1}$ and their union $J_\mathcal{Z}$ :
\begin{equation} \label {blocs}
J_k := \left\{ (s,y) \in \mathbb{R_+}\times \Z^2 , \ kn \leq s < (k+1)n \textrm{ and } |y-z_{k-1}\sqrt{n}| < C_1 \sqrt{n} \right\} \end{equation}
where $C_1$ will be a big constant.

\begin{center}
\begin{tikzpicture}[line cap=round,line join=round,>=triangle 45,x=2.0cm,y=0.8cm]
\draw [color=cqcqcq,dash pattern=on 1pt off 1pt, xstep=2.0cm,ystep=0.8cm] (0,-2.5) grid (5.93,4.88);
\clip(-0.65,-3) rectangle (5.93,4.88);
\fill[line width=0.4pt,color=qqffqq,fill=qqffqq,fill opacity=0.2] (0,2.5) -- (1,2.5) -- (1,-2.5) -- (0,-2.5) -- cycle;
\fill[line width=0.4pt,color=qqffqq,fill=qqffqq,fill opacity=0.2] (1,1.5) -- (2,1.5) -- (2,-3.5) -- (1,-3.5) -- cycle;
\fill[line width=0.4pt,color=qqffqq,fill=qqffqq,fill opacity=0.2] (2,4.5) -- (3,4.5) -- (3,-0.5) -- (2,-0.5) -- cycle;
\fill[line width=0.4pt,color=qqffqq,fill=qqffqq,fill opacity=0.2] (3,2.5) -- (4,2.5) -- (4,-2.5) -- (3,-2.5) -- cycle;
\fill[line width=0.4pt,color=qqffqq,fill=qqffqq,fill opacity=0.2] (4,5.5) -- (5,5.5) -- (5,0.5) -- (4,0.5) -- cycle;
\draw (0,0)-- (0.26,0.85)-- (0.8,-0.9)-- (1.29,0.64)-- (1.74,1.04)-- (2.1,2.78)-- (2.33,2.06)-- (2.59,0.94)-- (3.11,0.35)-- (3.4,2.69)-- (3.66,2.21)-- (4.08,3.7)-- (4.3,3.14)-- (5,0.6);
\draw [->] (0,0) -- (5.71,0);
\draw [->] (0.01,-1.58) -- (0,4.43);
\draw [line width=0.4pt,color=qqffqq] (0,2.5)-- (1,2.5);
\draw [line width=0.4pt,color=qqffqq] (1,2.5)-- (1,-2.5);
\draw [line width=0.4pt,color=qqffqq] (1,-2.5)-- (0,-2.5);
\draw [line width=0.4pt,color=qqffqq] (0,-2.5)-- (0,2.5);
\draw [line width=0.4pt,color=qqffqq] (1,1.5)-- (2,1.5);
\draw [line width=0.4pt,color=qqffqq] (2,1.5)-- (2,-3.5);
\draw [line width=0.4pt,color=qqffqq] (2,-3.5)-- (1,-3.5);
\draw [line width=0.4pt,color=qqffqq] (1,-3.5)-- (1,1.5);
\draw [line width=0.4pt,color=qqffqq] (2,4.5)-- (3,4.5);
\draw [line width=0.4pt,color=qqffqq] (3,4.5)-- (3,-0.5);
\draw [line width=0.4pt,color=qqffqq] (3,-0.5)-- (2,-0.5);
\draw [line width=0.4pt,color=qqffqq] (2,-0.5)-- (2,4.5);
\draw [line width=0.4pt,color=qqffqq] (3,2.5)-- (4,2.5);
\draw [line width=0.4pt,color=qqffqq] (4,2.5)-- (4,-2.5);
\draw [line width=0.4pt,color=qqffqq] (4,-2.5)-- (3,-2.5);
\draw [line width=0.4pt,color=qqffqq] (3,-2.5)-- (3,2.5);
\draw [line width=0.4pt,color=qqffqq] (4,5.5)-- (5,5.5);
\draw [line width=0.4pt,color=qqffqq] (5,5.5)-- (5,0.5);
\draw [line width=0.4pt,color=qqffqq] (5,0.5)-- (4,0.5);
\draw [line width=0.4pt,color=qqffqq] (4,0.5)-- (4,5.5);
\draw [line width=2pt,color=ffqqqq] (1,-1)-- (1,0);
\draw [line width=2pt,color=ffqqqq] (2,3)-- (2,2);
\draw [line width=2pt,color=ffqqqq] (3,1)-- (3,0);
\draw [line width=2pt,color=ffqqqq] (4,4)-- (4,3);
\draw [line width=2pt,color=ffqqqq] (5,0)-- (5,1);
\begin{scriptsize}
\draw[color=black] (1,0.22) node {$n$};
\draw[color=black] (3.02,-0.3) node {$3n$};
\draw[color=black] (5.05,-0.3) node {$5n$};
\draw[color=black] (-0.3,0.08) node {$0$};
\draw[color=black] (-0.3,1.06) node {$\sqrt{n}$};
\draw[color=black] (-0.3,2.05) node {$2\sqrt{n}$};
\draw[color=black] (-0.3,3.11) node {$3\sqrt{n}$};
\draw[color=black](-0.3,-0.93) node {$-\sqrt{n}$};
\draw[color=black] (2.07,-0.3) node {$2n$};
\draw[color=black] (4.07,-0.3) node {$4n$};
\draw[color=black] (5.87,0.15) node {$\R_+$};
\draw[color=black] (-0.3,4.14) node {$\Z$};
\end{scriptsize}
\end{tikzpicture}

the boxes $J_k$
\end{center}

Let $\delta>0$ be a small real number (we will fix its value in (\ref{delta})), we define a new measure $\wt{Q}_\z$ on the environment :
\[\frac{d\wt{Q}_\z}{dQ}= \prod_{k=0}^{m-1}\prod_{|y-z_k|<C_1\sqrt{n}} e^{-\delta (B_{(k+1)n}^y-B_{kn}^y)}=\prod_{y\in\Z}\exp\left( -\delta\int_0^T \ind{(s,y)\in J_\z} dB^y_s\right).\]
By Cameron-Martin-Grisanov theorem, under this new measure, the processes $(B^x)_{x\in\Z}$ are still independent and verify the following equations :
\[dB^x_t= d\wt{B}^x_t-\delta \ind{{(t,x)\in J_\z}} dt,\]
where the $(\wt{B}^x_t)_{x\in\Z}$ are independent Brownian Motions under $\wt{Q}_\z$.

We use H\"older's inequality to get an upper bound on $Q[\overline{Z}_\z^\theta]$ :
\begin{equation} \label{hold} \begin{array}{rcl}
Q[\overline{Z}_\mathcal{Z}^\theta] & = & \wt{Q}_\z\left[\left(\frac{d\wt{Q}_\z}{dQ}\right)^{-1}\overline{Z}_\mathcal{Z}^\theta\right] \\
& \leq & \wt{Q}_\z\left[\left(\frac{d\wt{Q}_\z}{dQ}\right)^{\frac{-1}{1-\theta}}\right]^{1-\theta} \wt{Q}_\z\left[ \overline{Z}_\mathcal{Z} \right]^\theta
\end{array} \end{equation}

we can compute the first term of (\ref{hold}), it is standard Gaussian calculus :
\[\wt{Q}_\z\left[\left(\frac{d\wt{Q}_\z}{dQ}\right)^{\frac{-1}{1-\theta}}\right]^{1-\theta} =Q\left[\left(\frac{d\wt{Q}_\z}{dQ}\right)^{\frac{-\theta}{1-\theta}}\right]^{1-\theta} =\exp\left(\frac{\delta^2\theta^2}{2(1-\theta)} |J_\z|\right).\]
We know that the size of $J_\z$ is approximately $2C_1 n\sqrt{n}m$. We will fix the value of $\delta$ such that this quantity is small (recall that $\theta=1/2$) :
\begin{equation} \label{delta} \textrm{let } \ \delta = C_1^{-1/2}n^{-3/4}, \ \textrm{ then } \wt{Q}_\z\left[\left(\frac{d\wt{Q}_\z}{dQ}\right)^{\frac{-1}{1-\theta}}\right]^{1-\theta}\leq e^{m}\end{equation}

\subsection{Upper bound on the second term}

Using (\ref{decomp2}) with (\ref{hold}) and (\ref{delta}), we obtain
\begin{equation} \label{decomp3} Q[Z_T^\theta] \leq e^m \sum_{\z \in \mathbb{Z}^m} \wt{Q}_\z [\overline{Z}_\z]^\theta .\end{equation}

Under the measure $\wt{Q}_\z$, all the brownian motion in the boxes $J_k$ have a negative drift $-\delta$, so \emph{if $X_t$ stays in the boxes}, $H_T(X)$ is a gaussian variable of mean  $-T\delta$ and variance $T$, and 
\[\wt{Q}_\z[\exp(\beta H_T(X)-T\beta^2/2)]=e^{-\beta \delta T}. \]
Now recall the value of $\delta$ from (\ref{delta}), we can take $n$ big enough such that this quantity is really small ($C_2$ will be a big constant):
\begin{equation}\label{n} \textrm{let } \ n = C_1^2 C_2\beta^{-4}, \ \textrm{ then }  e^{-\beta \delta n} =(e^{-\beta \delta T})^{1/m} \leq e^{-C_2^{1/4}} \ \textrm{ is really small.}\end{equation}

Now we have to consider that the trajectory may exit the box at some time. Let $\vep>0$. We note $T_{J_\z}=T_{J_\z}(X)$ the time spent in $J_\z$ by the trajectory $X$, and $T_{J_k}$ the time spent in the box $J_k$.
\ba  \wt{Q}_\z [\overline{Z}_\z]&=& \wt{Q}_\z P [\exp(\beta H_T(X)-T\beta^2/2)\ind{\{X_{in}\in I_{z_i}, i=1,\dots, m\} }]\\
&=& P[\exp(-\beta \delta T_{J_\z}) \ind{\{X_{in}\in I_{z_i}, i=1,\dots, m\} }]\\
&=& P\left[\prod_{k=0}^{m-1}\exp(-\beta \delta T_{J_k}) \ind{\{ X_{in}\in I_{z_i}, i=1,\dots, m\} }\right] \ea

A trajectory verifying $X_{in}\in I_{z_i}$ for $i=1,\dots ,m$ can be seen as $m$ slices of trajectory, respectively starting from one point of $I_{z_{i-1}}$ and finishing at a point of $I_{z_i}$, which are similar to slices starting from somewhere in $I_0$ and finishing in $I_{z_i-z_{i-1}}$. So we can bound the contribution in the above expectation by maximizing over the starting point $x$ in $I_0$, and we get the following upper bound 
\[ \wt{Q}_\z [\overline{Z}_\z]\leq \prod_{k=0}^{m-1} \max_{x\in I_0} P_x \left[\exp(-\beta \delta T_{J_0}) \ind{\{X_{n}\in I_{z_k-z_{k-1}}\} }\right]\]
where $P_x$ is the probability under which the trajectory starts from $x$.

Now that each term of the sum (\ref{decomp3}) is split in $m$ independent terms : 
\ba \sum_{\z \in \Z^m} \wt{Q}_\z [\overline{Z}_\z]^\theta &\leq & \sum_{\z \in \mathbb{Z}^m}\prod_{k=0}^{m-1} \max_{x\in I_0} P_x \left[\exp(-\beta \delta T_{J_0})\ind{\{ X_{n}\in I_{z_k-z_{k-1}}\} }\right]^\theta \\
& = &\left(\sum_{z\in\Z} \max_{x\in I_0} P_x \left[\exp(-\beta \delta T_{J_0}) \ind{\{ X_{n}\in I_z\} }\right]^\theta\right)^m.\ea

Basic properties of the random walk on $\Z$ tell us that the probability $P[X_n\in I_z]$ decreases like $\exp(-z^2/2)$ when $z\to\pm \infty$, so we can find a $R$ big enough such that 
\begin{equation} \label{eps1} \sum_{|z|>R} \max_{x\in I_0} P_x \left[\exp(-\beta \delta T_{J_0}) \ind{\{ X_{n}\in I_z\} }\right]^\theta <\vep.\end{equation}
Now we use the following trivial bound for the remainig terms :
\[\sum_{|z|\leq R} \max_{x\in I_0} P_x \left[\exp(-\beta \delta T_{J_0}) \ind{\{ X_{n}\in I_z\} }\right]^\theta\]
 \[\leq (2R+1)\max_{x\in I_0} P_x \left[\exp(-\beta \delta T_{J_0})\right]^\theta.\]

We now estimate roughly the time $T_{J_0}$ spent in the box : $T_{J_0}=n$ if the trajectory never leaves the box, and is bounded from below by 0 if it does :
\[P_x \left[\exp(-\beta \delta T_{J_0})\right] \leq \exp (-\beta \delta n)+ P_x[\textrm{the trajectory leaves the box}] \]

 Then we take $C_1$ big enough such that 
\begin{equation} \label{eps2} P\left[\max_{t\in[0,n]} |X_t|> (C_1-1)\sqrt{n}\right]<\frac{\vep}{2R+1},\end{equation}
and we just have to inject (\ref{eps1}) and (\ref{eps2}) in (\ref{decomp3}) and we get :
\[Q[Z_t^\theta] \leq e^m \left((2R+1)e^{-C_2^{1/4}}+2\vep\right)^m,\]
and taking $C_2$ big enough, this is smaller than $e^{-m}$. To conclude, we put this result in (\ref{fract}) and (\ref{hold}) :
\[\frac{1}{T} Q[\log Z_T]\leq \frac{-m}{\theta T}= \frac{-2}{n}=   \frac{-2 \beta^4}{C_1^2C_2}.\]

\section{Proof in dimension 2}

Once again we sketch the proof before starting :
\begin{enumerate}
\item We reduce the problem to the exponential decay of $Q[Z_t^\theta]$ for some $\theta \in (0,1)$.
\item We use a decomposition of $Z_t$ by splitting it into different contributions that corresponds to trajectories that stays in a large corridor.
\item To estimate the fractional moment terms appearing in the decomposition, we will add a term that will handicap the environments that have a lot of positive correlations around the corridors corresponding to each contribution. We define this handicapping term in such a way that the new measure is not very different from the original one.
\item We use some basic properties of a Random Walk in $\Z^2$ to compute the expectation with the new term.
\end{enumerate}

The steps 1 and 2 are exactly the same than in the 1-dimensional case, refer to the sections 2.1 and 2.2. We define the interval $I_z$ as follows :
\[\textrm{if }\ z=(a,b)\in \Z^2, \ \textrm{ then }\ I_z=[a\sqrt{n},(a+1)\sqrt{n})\times [b\sqrt{n}, (b+1)\sqrt{n}).\]
We obtain the following bound :
\begin{equation} \label{decomp2'} 
Q[Z_T^\theta] \leq \sum_{\z \in (\mathbb{Z}^2)^m} Q[\overline{Z}_\mathcal{Z}^\theta ] .\end{equation}

\subsection{The tweaking of the measure}

We fix one $\mathcal{Z}=(z_1,\dots,z_m)$. We define $m$ blocs $(J_k)_{k=0}^{m-1}$ and their union $J_\mathcal{Z}$ :
\begin{equation} \label {blocs'}
J_k := \left\{ (s,y) \in \mathbb{R_+}\times \Z^2 , \ kn \leq s < (k+1)n \textrm{ and } |y-z_{k-1}\sqrt{n}| < C_3 \sqrt{n} \right\} \end{equation}
where $C_3$ will be a big constant. In the following, I will abuse this notation and say that $y\in J_k$ if $|y-z_k\sqrt{n}|<C_3\sqrt{n}$ (\emph{ie} $y$ is in the vertical component of $J_k$).

First of all, let me go on a little tangent. If $B$ and $\widehat{B}$ are two brownian motions, let us define what I will call \emph{the correlation between $B$ and $\widehat{B}$ in $[0,T]$} $\rho_T(B,\widehat{B})$ :
\[\rho_T(B,\widehat{B}) =\int_0^T \left(\int_0^s dB_u\right) d\widehat{B}_s + \int_0^T \left(\int_0^s d\widehat{B}_u\right) dB_s .\]
This quantity is easier to understand if we consider that the Brownian motions are constituted of a sum of tiny discrete increments $\bigtriangleup_{n} B_s=B_{s+\frac{1}{n}}-B_s$. Then the quantity becomes 
\[\rho_T(B,\widehat{B}) = \lim_{n\to \infty} \sum_{\begin{array}{c} 0\leq s,u<T , s\neq u\\ s,u\in \frac{1}{n}\Z \end{array}} \bigtriangleup_{n} B_s \bigtriangleup_{n} \widehat{B}_u.\]
So if the quantity $\rho_T(B,\widehat{B})$ is big, it means that there are a lot of discrete increments  $(\bigtriangleup_{n} B_s, \bigtriangleup_{n} \widehat{B}_u)$ that are the same sign, or in rough term that $B$ and $\widehat{B}$ are both ``mostly increasing'' or both ``mostly decreasing'' (actually, it is easy to compute that $\rho_T(B,\wh{B})=B_T\cdot \wh{B}_T$, but I wanted to stress the ``tiny increments'' vision).

The quantity that will measure the correlations of the environment in one box $J_k$ is a more complicated version of what we just did :
\begin{equation} \label{cordef'}
R_k =R_k(\omega) = \sum_{x\in J_k} \sum_{y\in J_k} \int_{kn}^{(k+1)n} \left(\int_{kn}^s  V_{(s,x)(u,y)}  dB^y_u\right) dB^x_s  \end{equation}
where
\begin{equation} \label{vv'}
V_{(s,x)(u,y)} := \frac{\ind{ \{ |x-y| \leq C_4 \sqrt{|s-u|} \}}}{100 C_3 C_4 n \sqrt{\log n} (|s-u|+ 1)} . \end{equation}

The quantity $R_k$ has mean zero, and, roughly speaking, it is positive when there is a lot of close brownian motions which are ``mostly going the same way". To compute the mean and variance of $R_k$ we simply use that the $B^x_t$ are independent Brownian motions, on the same filtration $\g_t$, and then
\[Q\left[\left(\int_0^n f(s) dB^x_s\right) \left(\int_0^n g(u) dB^y_u \right) \right] = \int_0^n f(s)g(s) ds \ \textrm{ if }\ x=y \ \textrm{ and }\ 0 \ \textrm{ otherwise}.\]
\[Var_Q(R_k) = \sum_{x\in J_k} \sum_{y\in J_k} \int_{kn}^{(k+1)n} \int_{kn}^s  V_{(s,x)(u,y)}^2  du ds.\]
Let us compute this term :
\ba Var_Q(R_0) & = &
\sum_{x\in J_0} \sum_{y\in J_0} \int_0^n \int_0^s \frac{\ind{ \{ |x-y| \leq C_4 \sqrt{|s-u|} \}}}{10^4 C_3^2 C_4^2 n^2 \log n (|s-u|+ 1)^2}  du ds  \\
&\leq & \sum_{x\in J_0}  \int_0^n \int_0^s  \frac{4C_4^2|s-u|}{10^4 C_3^2 C_4^2 n^2 \log n (|s-u|+ 1)^2}  du ds \\ 
&\leq & 4C_3^2n \int_0^n \int_0^s  \frac{4C_4^2|s-u|}{10^4 C_3^2 C_4^2 n^2 \log n (|s-u|+ 1)^2}  du ds \\
&\leq  &4C_431^2 n \int_0^n  \frac{4C_4^2}{10^4 C_3^2 C_4^2 n^2 \log n }  \log(n) ds\\% 
&\leq &\frac{16  C_3^2 C_4^2 n^2 \log n}{10^4 C_3^2 C_4^2 n^2 \log n}\leq 1  \ea

Let $K$ be a large constant. One defines the function $f_K$ on $\mathbb{R}$ to be :
\[f_K(x) := -K \ind{ \{x > \exp K^2\} }. \]
We define $g_\mathcal{Z}$ function of the environment as
\[ g_\mathcal{Z}(\omega) := \exp \left( \sum_k f_K\left(R_k \right) \right) \]
So if in a box $J_k$ the quantity $R_k$ is too big, we affect a weight $e^{-K}$. 

We use H\"older's inequality to get an upper bound on $Q[\overline{Z}_\z^\theta]$ :
\begin{equation} \label{hold'} \begin{array}{rcl}
Q[\overline{Z}_\mathcal{Z}^\theta] & = & Q[g_\mathcal{Z}(\omega)^{-\theta}(g_\mathcal{Z}(\omega)\overline{Z}_\mathcal{Z})^\theta] \\
& \leq & Q[g_\mathcal{Z}(\omega)^{- \frac{\theta}{1-\theta}}]^{1-\theta} Q[g_\mathcal{Z}(\omega) \overline{Z}_\mathcal{Z} ]^\theta
\end{array} \end{equation}

The block structure of $g_\mathcal{Z}$ allows to express the first term as a power of $m$.
\begin{equation} 
 Q\left[ g_\mathcal{Z}(\omega)^{-\frac{\theta}{1-\theta}} \right] = Q \left[ \exp \left( - \frac{\theta}{1-\theta} f_K\left(  R_0 \right) \right) \right]^m
\end{equation}

We have already seen that $Q[R_0]=0$ and $Var_Q(R_0)\leq 1$. In consequence,
\[ Q\left[ R_0  > \exp K^2 \right] \leq  \exp (-2K^2), \]
hence 
\begin{equation} \label{firstterm'}
Q\left[g_\mathcal{Z}(\omega)^{- \frac{\theta}{1-\theta}}\right] \leq \left(1+\exp\left(\frac{\theta}{1-\theta} K-2K^2\right)\right)^m \leq 2^m
\end{equation}
if $K$ is large enough. We are left with estimating the second term
\begin{equation} \label{secondterm'}
Q\left[ g_\mathcal{Z}(\omega) \overline{Z}_\mathcal{Z} \right] = Q P \left[g_\mathcal{Z}(\omega)\exp (\beta H_T-T\beta^2/2) \ind{\{X_{in} \in I_{z_i}\ ,\ i=1 \dots m \} }\right] \end{equation}

For a fixed trajectory of the random walk $X$, we consider $\widehat{Q}_X$ the modified measure on the environment with density 
\begin{equation} \frac{\textrm{d}\widehat{Q}_X}{\textrm{d}Q} := \exp (\beta H_T(X)-T\beta^2/2) . \end{equation}

Under the measure $\widehat{Q}_X$, the $B^x_t$ are still independent processes on the filtration $\g_t$, and 
\[ dB^x_t= d \widehat{B}^x_t + \beta \ind{\{X_t=x\} } dt,\]
where the $\widehat{B}^x_t$ are independent brownian motions under $\widehat{Q}_X$.

With the measure $\widehat{Q}_X$, (\ref{secondterm'}) becomes
\[Q\left[ g_\mathcal{Z}(\omega) \overline{Z}_\mathcal{Z} \right] =  P \widehat{Q}_X [g_\mathcal{Z}(\omega) \ind{ \{ X_{in} \in I_{z_i}  , \ i=1 \dots m\} }]. \]
Now we try to bound this term by something that could be decomposed on the blocks $J_k$.

As we did in dimension 1, we slice the trajectory $X_{[0,T]}$ in $m$ bits of trajectories starting from somewhere in $I_0$ and finishing in $I_{z_i-z_{i-1}}$. So we can bound the contribution in the above expectation by maximizing over the starting point $x$ in $I_0$, and we get the following upper bound 

\begin{equation} \label{upper'} 
Q\left[ g_\mathcal{Z}(\omega) \overline{Z}_\mathcal{Z} \right] \leq \prod_{i=1}^m \max_{x \in I_0} P_x \widehat{Q}_X \left[\exp \left(f_K \left( R_0 \right) \right) \ind{ \{ X_n \in I_{z_i-z_{i-1}} \} }\right]. \end{equation}

Using this with (\ref{firstterm'}),  (\ref{hold'}) and (\ref{decomp2'}) we get the inequality
\[ P[Z_t^\theta] \leq 2^{m(1-\theta)}\left( \sum_{z \in \mathbb{Z}^2}  \max_{x \in I_0} P_x \widehat{Q}_\omega \left[\exp \left(f_K \left( R_0 \right) \right) \ind{ \{ X_n \in I_z \} } \right]^\theta \right) ^m.\]
Therefore to prove the exponential decay of $Q[W_t^\theta]$ it is sufficient to prove that
\[ \sum_{z \in \mathbb{Z}^2} \max_{x \in I_0} P_x \widehat{Q}_X \left[ \exp \left(f_K \left( R_0 \right) \right) \ind{ \{ X_n \in I_z \} } \right]^\theta \]
is small.

\subsection{Computation of the second term}

We fix some $\vep > 0$. Asymptotic properties of the random walk guarantee that the probability $Q[X_n \in I_z]$ decays like $e^{-|z|^2/2}$ when $|z| \to \infty$. So we can find $R$ such that
\begin{equation} \label{far'} \sum_{|z|>R}  \max_{x \in I_0} P_x \widehat{Q}_X \left[ \exp \left(f_K \left( R_0\right) \right) \ind{ \{ X_n \in I_z \} } \right]^\theta  
\leq \sum_{|z|> R} \max_{x \in I_0} P_x [X_n \in I_z]^\theta \leq \epsilon \end{equation}
To estimate the rest of the sum, we have the trivial bound
\[ \sum_{|z|\leq R} \left[ \max_{x \in I_0} P_x \widehat{Q}_X \left[\exp \left(f_K \left(R_0 \right) \right)\right] \ind{ \{ X_n \in I_{z_i-z_{i-1}} \} } \right]^\theta \]
\begin{equation} \label{close'} \leq 4R^2 \max_{x \in I_0} P_x \widehat{Q}_X \left[ \exp \left(f_K \left( R_0 \right) \right)\right] ^\theta \end{equation}
Therefore it is sufficient for our purpose to check 
\begin{equation} \label{toto'} \max_{x \in I_0} P_x \widehat{Q}_X \left[ \exp \left(f_K \left( R_0 \right) \right)\right] < \vep' \end{equation}
for some small $\vep'$.

Now we have to find an upper bound for
\[ P_x \widehat{Q}_X \left[ \exp \left(f_K(R_0) \right)\right] \]
that is uniform in $x \in I_0$. To study the left-hand term of the equation (\ref{firstterm'}), we showed that under $Q$, $R_0$ had a mean 0 and a small variance, so $f_K(R_0)$ was 0 with very high $Q$-probability. Now we will show that under $\widehat{Q}_X$, $R_0$ has a really high mean and a small variance, so $f_K(R_0)=-K$ with high $\widehat{Q}_X$-probability.

We fix $\delta>0$ small. First, we force the random walk to stay in the zone where the environment is modified by writing
\[ \max_{x\in I_0} P_x \widehat{Q}_X \left[ \exp \left(f_K \left( R_0 \right) \right)\right] \leq \]
\begin{equation} \label{kali'} 
P_0\left[\max_{s \in [0,n]} |X_s| \geq (C_3-1) \sqrt{n} \right]
+\max_{x\in I_0} P_x \widehat{Q}_X \left[ \exp \left(f_K \left( R_0 \right) \right)\ \ind{ \{ X_{[0,n]}\subseteq J_0\} }\ \right]\end{equation}
and we can take $C_3$ big enough such that $ P_0[\max_{[0,n]} |X_s| \geq (C_3-1) \sqrt{n} ] \leq \delta$, so we can concentrate on the trajectories that are inside the box $J_0$.

Recall that under $\widehat{Q}_X$, $dB^x_t= d \widehat{B}^x_t + \beta \ind{{X_t=x}} dt$. Now we can compute the expectation of $R_0$:
\begin{equation} \label{exp'}
 \widehat{Q}_X \left[ R_0 \right] = \beta^2 \int_0^n\int_0^s V_{(s,X_s)(u,X_u)}du \ ds  \end{equation}

To study this term, let us define the quantity
\begin{equation} \label{dndef'} 
D_n:= \int_0^n\int_0^s \frac{1}{n \sqrt{\log n} (|s-u|+ 1)}du\ ds \end{equation}
and the random variable (under $P$) $Y$ :
\begin{equation} \label{xdef'}
Y:= \int_0^n\int_0^s  \frac{\ind{ \{ |X_s-X_u|\leq C_4 \sqrt{|s-u|}\}} }{n \sqrt{\log n} (|s-u|+ 1)}du\ ds \end{equation}
Roughly speaking, $D_n$ is what $Y$ would be if $C_4 = \infty$. It is easy to check, with the central limit theorem, that we can take $C_4$ large enough such that $P[Y] \geq (1-\delta)D_n$, uniformly in $n$.

Then, as $Y \leq D_n$ a.s. it implies that 
\begin{equation} \label{mark'}
P_x \left[ Y \leq \frac{D_n}{2} \right] \leq 2 \delta \end{equation}
and  a quick computation of the value of $D_n$ gives us :
\[ D_n \geq \sqrt{\log n} . \]
If we inject this result in the term in (\ref{exp'}), we have that with large $P_x$-probability,
\begin{equation} \label{exp2'}
\widehat{Q}_X \left[ R_0 \right] \geq \frac{\beta^2 \sqrt{\log n} }{200 C_3 C_4} \end{equation}
Now we can use the only free hand we have left : the value of $n$ 
\begin{equation} \label{n'}\textrm{let }\ n= \exp\left(\frac{C_5}{\beta^4}\right), \ \textrm{ then }\ \widehat{Q}_X \left[ R_0 \right] \geq \frac{\sqrt{C_5} }{200 C_3 C_4}\geq 2\exp K^2\end{equation}
if we take $C_5$ large enough.

Now we have to  bound the variance of $R_0$ under $\widehat{Q}_X$. We can decompose $R_0$, 
\[\sum_{x,y\in J_0} \int_0^n \left(\int_0^s  V_{(s,x)(u,y)}  dB^y_u\right) dB^x_s
 = \beta^2 \int_0^n\int_0^s V_{(s,X_s)(u,X_u)} du\ ds \]
\[ + \beta\sum_{x\in J_0} \int_0^n \int_0^n V_{(s,x)(u,X_u)}du\ d\widehat{B}^x_s   
 +\sum_{x,y\in J_0} \int_0^n \left(\int_0^s  V_{(s,x)(u,y)}  d\widehat{B}^y_u\right) d\widehat{B}^x_s\]
Now, using that $(x+y)^2 \leq 2x^2+2y^2$, we can bound the variance by

\[2\beta^2 n \sum_{x\in J_0}\left(\int_0^n V_{(s,x)(u,X_u)}du\right)^2 + 2 \sum_{x\in J_0}\sum_{y\in J_0} \int_0^n \int_0^s  V_{(s,x)(u,y)}^2  du\ ds.\]
We already computed the right-hand term when we studied $Var_Q(R_0)$, it is $\leq 1$. With similar techniques we find the other one :
\ba \sum_{x\in J_0}\int_0^n V_{(s,x)(u,X_u)}du &\leq & \frac{C_4}{10C_3\sqrt{\log n}}, \\
\sup_{x\in J_0}\int_0^n V_{(s,x)(u,X_u)}du &\leq & \frac{\sqrt{\log n}}{10C_3C_4n}, \\
\textrm{therefore } \ 2\beta^2 n \sum_{x\in J_0}\left(\int_0^n V_{(s,x)(u,X_u)}du\right)^2&\leq & \frac{2\beta^2}{100 C_3^2} \leq 1\ea

Now, using (\ref{kali'}), (\ref{mark'}) and (\ref{n'}), we can compute (\ref{toto'}) :
\begin{equation} P_x \widehat{Q}_X\left| \exp (f_K( R_0)) \leq P_0\left[\max_{s \in [0,n]} |X_s| \geq (C_3-1) \sqrt{n} \right] \right]\]

\[  P_x\left[ \widehat{Q}_X[R_0] \leq 2e^{ K^2}\ind{ \{ X_{[0,n]}\subseteq J_0\} }\right] + P_x\widehat{Q}_X \left[\exp (f_K(R_0)) \ind{\{ \widehat{Q}_X[R_0] > 2e^{ K^2}\} }\right] \end{equation}

\[\leq 2\delta + e^{-K} + 3 e^{- 2K^2} \]
so that our result is proved provided that $K$ and $\delta$ have been chosen respectively large and small enough. 

In conclusion, we have proved (\ref{toto'}), which with  (\ref{decomp2'}), (\ref{far'}) and (\ref{close'}) allow us to say 
\[Q[Z_T^\theta]\leq e^{-m}\]
and with (\ref{fract}) and the value of $n$ from (\ref{n'}) :
\[\frac{1}{T} Q[\log Z_T]\leq \frac{-m}{\theta T}= \frac{-2}{n}= \exp\left( \frac{-c}{\beta^4}\right).\]

\subsection*{Aknowledgements}

I would like to thank my PhD advisor Francis Comets for his careful reading and numerous suggestions, and without whom this paper would not have been possible.


\begin{thebibliography}{99}
\bibitem{cativi}{A.~Cadel, S.~Tindel, F.~Viens : Sharp asymptotics for the partition function of some continuous-time directed polymers. \emph{Potential Anal.} 29 (2008) 139-166}
\bibitem{cahu}{P.~Carmona, Y.~Hu : Strong disorder implies strong localization dor directed polymers in a random environment, \emph{ALEA Lat. Am. Math. Stat.} 2 (2006) 217-229}
\bibitem{cako}{P.~Carmona, L.~Koralov, S.~Molchanov : Asymptotics for the almost sure Lyapunov exponent for the solution of the parabolic Anderson model, \emph{Rand. Oper. Sto. equ.} 9 (2001) no1, 77-86}
\bibitem{camo}{P.~Carmona, S.~Molchanov : Parabolic Anderson problem and intermittency, \emph{Mem. Amer. Math. Soc.} 108 (1994)}
\bibitem{cocr}{F.~Comets, M.~Cranston : Overlaps and pathwise localization in phe Anderson polymer model, preprint 2011}
\bibitem{coyo}{F.~Comets, N.~Yoshida : Brownian directed polymers in random environment, \emph{Comm. Math. Phys.} $\mathbf{254}$ (2005), 257-287}
\bibitem{crgamo}{M.~Cranston, D.~Gauthier, T.~Mountford : On large deviations for the parabolic Anderson model, \emph{Proba. Theo. Related Fields} 147 (2010) 349-378}
\bibitem{crmosh}{M.~Cranston, T.~Mountford, T.~Shiga : Lyapunov exponents for the parabolic Anderson model, \emph{Acta Math. Univ. Comenian. (N.S.)} 71 (2002) 163-188}
\bibitem{gaho}{J.~G\"artner, F.~den Hollander : Intermittency in a catalytic random medium, \emph{Ann. Probab.} 34 (2006) 2219-2287}
\bibitem{lac}{H.~Lacoin : New bounds for the free energy of directed polymers in dimension 1+1 and 1+2,  \emph{Comm. Math. Phy.} 294 (2010) 471-503}
\bibitem{lac2}{H.~Lacoin : Non-coincidence of Quenched and Annealed Connective Constants on the supercritical planar percolation cluster, arXiv: 1203.6051}
\bibitem{mo}{G.~Moreno : Asymmetric directed polymers in random environments, arXiv: 1009.5576}
\bibitem{mooc}{J.~Moriarty, N.~O'Connell : On the free energy of a directed polymer in a Brownian environment, \emph{Markov Proc. Rel. Fields} 13 (2007) 251-266}
\bibitem{yo1}{N.~Yoshida : Phase transitions for the growth rate of Linear Stochastic Evolutions, preprint 2008, \emph{J. Stat. Phys.} 133 No6 1033-1058}
\bibitem{yo2}{N.~Yoshida : Localization for Linear Stochastic Evolutions, \emph{J. Stat. Phys.} 138 , No.4/5, 598--618, (2010).}

\end{thebibliography}
\end{document}